\documentclass[11pt]{article}
\usepackage{epsfig}
\usepackage{amssymb,amsmath,amsthm,amscd}
\usepackage{latexsym}

\newtheorem{thm}{Theorem}[section]
\newtheorem{prop}[thm]{Proposition}
\newtheorem{cor}[thm]{Corollary}
\newtheorem{lem}[thm]{Lemma}

\theoremstyle{definition}
\newtheorem{defn}[thm]{Definition}

\theoremstyle{remark}

\newcounter{labelflag} 

\newcommand{\Label}[1]{
                       \ifnum\thelabelflag=1
                          \ifmmode
                             \makebox[0in][l]{\qquad\fbox{\rm#1}}
                          \else
                             \marginpar{\vspace{0.7\baselineskip}
                                        \hspace{-1.1\textwidth}
                                        \fbox{\rm#1}}
                          \fi
                       \fi
                       \label{#1} }

\newcommand{\be}{\begin{equation}}
\newcommand{\ee}{\end{equation}}

 \newcommand{\R}{\mathbb{R}}

 \def  \ltwo {L^2 (\R^n)}

\begin{document}

\begin{titlepage}
\title{\Large\bf  Pullback Attractors for  Non-autonomous
Reaction-Diffusion Equations  on  $\R^n$}
\vspace{7mm}

\author{
Bixiang Wang  \thanks {Supported in part by NSF  grant DMS-0703521.}
\vspace{5mm}\\
Department of Mathematics, New Mexico Institute of Mining and
Technology \vspace{1mm}\\ Socorro,  NM~87801, USA \vspace{5mm}\\
Email: bwang@nmt.edu \qquad Fax: (1-505) 835 5366}
\date{}
\end{titlepage}

\maketitle

\medskip

\begin{abstract}
We study  the long time behavior of solutions  of the
non-autonomous  Reaction-Diffusion  equation
defined on the entire space $\R^n$ when external terms
are unbounded in  a phase space.  The existence of a
pullback global attractor  for the equation  is established
in $L^2(\R^n)$ and $H^1(\R^n)$, respectively.
The pullback asymptotic
compactness of  solutions   is proved
 by using uniform a priori estimates
 on the tails of solutions   outside  bounded domains.
\end{abstract}



\baselineskip=1.4\baselineskip

\section{Introduction}
\setcounter{equation}{0}

In this paper,  we study the dynamical behavior of the
 non-autonomous  Reaction-Diffusion equation  defined on $\R^n$:
\be
  \label{intr1}
 \frac {\partial u}{\partial t}
 -   \Delta u + \lambda  u  = f(x, u) + g(x,t),
 \ee
 where $\lambda$   is a  positive constant,
 $g$  is  a  given function  in $L^2_{loc} (\R, L^2(\R^n))$,  and $f$ is a
nonlinear
function  satisfying a dissipative condition.

 Global attractors for non-autonomous  dynamical systems
 have been extensively studied in the literature, see, e.g.,  \cite{ant,
aul1, car1, car2, car3, car4,
cheb1, cheb2, cheb3, che, har,
lan1, lan2,  lu,  moi, pri2,
sun1, sun2, wan2, wangy}.
Particularly,  when  PDEs
 are defined in bounded domains, such attractors
 have been investigated in
 \cite{car1, car4, cheb1, cheb2, che, har, lu, sun1, sun2}.
 In the case of unbounded domains,  global attractors
 for non-autonomous  PDEs
 have been examined in \cite{ant,moi, pri2}
 for   almost periodic external terms,
 and in
     \cite{car2, car3, wan2, wangy} for unbounded external terms.
 In this paper, we  will prove   existence
 of a pullback   attractor for   equation
 \eqref{intr1}
  defined on   $\R^n$ with unbounded external terms.

Notice that the domain $\R^n$ for    \eqref{intr1}   is    unbounded,
and hence Sobolev
 embeddings are no longer compact in this case.
 This   introduces a major  obstacle
 for examining  the    asymptotic compactness of solutions.
 For  some PDEs,  such  difficulty     can be overcome
by  the  energy equation approach,  which was  introduced  by Ball in
\cite{bal1, bal2}  (see  also
  \cite{ car2, car3, gou1, ju1, luk,  moi, moi2,  ros1, wan5, wanx}).
   In this paper, we will  use the
  uniform estimates on the  tails  of solutions to
 circumvent the difficulty caused by the unboundedness of the domain.
   This idea was developed in \cite{wan}
 to prove asymptotic compactness of  solutions  for autonomous
parabolic  equations on $\R^n$,
and   later  extended to non-autonomous  equations
 in \cite{ant, pri2, wan2, wangy} and stochastic equations in \cite{bat3, wan3, wan4}.
 Here, we will use the method of tail-estimates
   to  investigate  the asymptotic behavior
 of equation \eqref{intr1} with nonlinearity of arbitrary   growth rate.
 We first  establish    the  pullback asymptotic
compactness of  solutions  of  equation \eqref{intr1}   and  prove   existence
of a pullback  global attractor  in $L^2(\R^n)$.
 Then we extend this result and   show
   existence of a pullback global attractor
   in $H^1(\R^n)$.

   It is worth noticing that   attractors
   for the non-autonomous Reaction-Diffusion equation defined
   on $\R^n$ with unbounded external terms were
   also studied in
   \cite{wangy}, where the authors proved the existence of
   a pullback  attractor when the nonlinearity $f$ satisfies  a Sobolev growth
   rate. In the present  paper, we deal
    with   the case where  the growth order of  $f$
    is arbitrary.
    The asymptotic compactness of solutions in \cite{wangy}
    was obtained by using the energy equation approach. But, here we will
    derive such compactness directly from the
     uniform tail-estimates of solutions. As we will see later,
     the existence of an attractor in $H^1(\R^n)$ is an immediate
     consequence of the existence of an attractor
     in $L^2(\R^n)$ and the asymptotic compactness of solutions
      in $H^1(\R^n)$.

 The paper is organized  as follows. In the next section, we
 recall   fundamental concepts and  results
for  pullback
   attractors  for non-autonomous  dynamical systems.
   In Section 3,   we   define a cocycle  for the non-autonomous
 Reaction-Diffusion
    equation  on $\R^n$.
   Section 4  is devoted  to deriving  uniform estimates of solutions
    for large space and time
   variables. In the last section, we
prove the   existence of a pullback
    global attractor for the equation in $L^2(\R^n)$ and
$H^1(\R^n)$.

   The following notations will be used throughout the paper.
 We denote by
$\| \cdot \|$ and $(\cdot, \cdot)$ the norm and   inner product
in  $L^2(\R^n)$ and use $\| \cdot\|_{p}$    to denote   the norm  in
$L^{p}(\R^n)$.    Otherwise, the
norm of  a general  Banach space $X$  is written as    $\|\cdot\|_{X}$.
 The letters $C$ and $C_i$ ($i=1, 2, \ldots$)
are  generic positive constants  which may change their  values from line to
line or even in the same line.

\section{Preliminaries}
\setcounter{equation}{0}

In this section,  we recall some basic concepts
related to pullback attractors    for non-autonomous  dynamical
systems.  It is  worth  noticing  that these concepts
are quite similar to that of random attractor  for stochastic
systems.  We refer the reader to \cite{arn1, bat1, bat3,  car22,  car2, car3,  cheb1, chu,   fla1, sun1, wan3}
for more details.

Let $\Omega$  be   a nonempty set and  $X$   a metric space with  distance
$d(\cdot, \cdot)$.

\begin{defn}
A   family of mappings $\{\theta_t\}_{t\in \R}$
from $\Omega$ to itself is called a  family of shift operators on $\Omega$ if
 $\{\theta_t\}_{t\in \R}$ satisfies the  group properties:

   (i) \  $\theta_0 \omega =\omega,  \quad  \forall \ \omega \in \Omega;$

   (ii)\ $  \theta_t (\theta_\tau \omega) = \theta_{t+\tau}  \omega,  \quad
  \forall \ \omega \in \Omega  \quad \mbox{and} \ \  t, \ \tau \in \R.$
  \end{defn}

\begin{defn}
Let $\{\theta_t\}_{t\in \R}$
be a   family of shift operators on $\Omega$.  Then a  continuous $\theta$-cocycle
$\phi$ on $X$
is  a mapping
$$
\phi: \R^+ \times \Omega \times X \to X, \quad (t, \omega, x) \mapsto \phi(t, \omega, x),
$$
which  satisfies, for  all  $\omega \in \Omega$ and
$t,  \tau \in \R^+$,

(i) \  $\phi(0, \omega, \cdot) $ is the identity on $X$;

(ii) \  $\phi(t+\tau, \omega, \cdot) = \phi(t, \theta_\tau \omega, \cdot) \circ \phi(\tau, \omega, \cdot)$;

(iii) \  $\phi(t, \omega, \cdot): X \to  X$ is continuous.
\end{defn}

Hereafter, we always assume that
$\phi$ is a continuous $\theta$-cocycle on $X$, and   $\mathcal{D}$ a  collection  of families of subsets of $X$:
$$
{\mathcal{D}} = \{ D =\{D(\omega)\}_{\omega \in \Omega}: \ D(\omega) \subseteq X
\  \mbox{for every} \ \omega \in \Omega \}.
$$
Such a collection ${\mathcal{D}}$ is often referred to as a universe in the literature.

\begin{defn}
Let $\mathcal{D}$ be a collection of families of  subsets of $X$.
Then  $\mathcal{D}$ is called inclusion-closed if
   $D=\{D(\omega)\}_{\omega \in \Omega} \in {\mathcal{D}}$
and  $\tilde{D}=\{\tilde{D}(\omega) \subseteq X:  \omega \in \Omega\} $
with
  $\tilde{D}(\omega) \subseteq D(\omega)$ for all $\omega \in \Omega$ imply
  that  $\tilde{D} \in {\mathcal{D}}$.
  \end{defn}

\begin{defn}
Let $\mathcal{D}$ be a collection of families of  subsets of $X$ and
$\{K(\omega)\}_{\omega \in \Omega} \in \mathcal{D}$. Then
$\{K(\omega)\}_{\omega \in \Omega} $ is called a  pullback
 absorbing
set for   $\phi$ in $\mathcal{D}$ if for every $B \in \mathcal{D}$
and   $\omega \in \Omega$, there exists $t(\omega, B)>0$ such
that
$$
\phi(t, \theta_{-t} \omega, B(\theta_{-t} \omega)) \subseteq K(\omega)
\quad \mbox{for all} \ t \ge t(\omega, B).
$$
\end{defn}

\begin{defn}
 Let $\mathcal{D}$ be a collection of families of  subsets of $X$.
 Then
$\phi$ is said to be  $\mathcal{D}$-pullback asymptotically
compact in $X$ if  for  every  $\omega \in \Omega$,
$\{\phi(t_n, \theta_{-t_n} \omega,
x_n)\}_{n=1}^\infty$ has a convergent  subsequence  in $X$
whenever
  $t_n \to \infty$, and $ x_n\in   B(\theta_{-t_n}\omega)$   with
$\{B(\omega)\}_{\omega \in \Omega} \in \mathcal{D}$.
\end{defn}

\begin{defn}
 Let $\mathcal{D}$ be a collection of families of  subsets of $X$
 and
 $\{\mathcal{A}(\omega)\}_{\omega \in \Omega} \in {\mathcal{D}}$.
Then     $\{\mathcal{A}(\omega)\}_{\omega \in \Omega}$
is called a    $\mathcal{D}$-pullback global  attractor  for
  $\phi$
if the following  conditions are satisfied,  for every  $\omega \in \Omega$,

(i) \  $\mathcal{A}(\omega)$ is compact;

(ii) \ $\{\mathcal{A}(\omega)\}_{\omega \in \Omega}$ is invariant, that is,
$$ \phi(t, \omega, \mathcal{A}(\omega)  )
= \mathcal{A}(\theta_t \omega), \ \  \forall \   t \ge 0;
$$

(iii) \ \ $\{\mathcal{A}(\omega)\}_{\omega \in \Omega}$
attracts  every  set  in $\mathcal{D}$,  that is, for every
 $B = \{B(\omega)\}_{\omega \in \Omega} \in \mathcal{D}$,
$$ \lim_{t \to  \infty} d (\phi(t, \theta_{-t}\omega, B(\theta_{-t}\omega)), \mathcal{A}(\omega))=0,
$$
where $d$ is the Hausdorff semi-metric given by
$d(Y,Z) =
  \sup_{y \in Y }
\inf_{z\in  Z}  \| y-z\|_{X}
 $ for any $Y\subseteq X$ and $Z \subseteq X$.
\end{defn}

The following existence result  of  a    pullback global  attractor
for a  continuous cocycle
can be found in \cite{arn1, bat1, bat3,  car22,  car2, car3,  cheb1, chu,   fla1}.

\begin{prop}
\label{att} Let $\mathcal{D}$ be an  inclusion-closed  collection of families of   subsets of
$X$ and $\phi$ a continuous $\theta$-cocycle on $X$.
 Suppose  that $\{K(\omega)\}_{\omega
\in \Omega} \in {\mathcal{D}} $ is a   closed   absorbing set  for  $\phi$  in
$\mathcal{D}$ and $\phi$ is $\mathcal{D}$-pullback asymptotically
compact in $X$. Then $\phi$ has a unique $\mathcal{D}$-pullback global
attractor $\{\mathcal{A}(\omega)\}_{\omega \in \Omega} \in {\mathcal{D}}$ which is
given by
$$\mathcal{A}(\omega) =  \bigcap_{\tau \ge 0} \  \overline{ \bigcup_{t \ge \tau} \phi(t, \theta_{-t} \omega, K(\theta_{-t} \omega)) }.
$$
\end{prop}

\section{Cocycle associated  with the  Reaction-Diffusion equation}
\setcounter{equation}{0}

In this  section,  we   construct a $\theta$-cocycle $\phi$
for the
 non-autonomous  Reaction-Diffusion equation
   defined on $\R^n$.
  For every $\tau \in\R$ and $t > \tau$, consider  the problem:
\be
\label{rd1}
    \frac{\partial u}{\partial t} - \Delta  u + \lambda u = f(x,u) + g(x, t),  \quad  x \in
\R^n ,
\ee
with the initial condition
\be
\label{rd2}
    u(x,\tau) = u_{\tau}(x), \hspace{3 mm}  x \in \R^n,
\ee
where $\lambda$ is a positive constant, $g $ is
given in $ L^2_{loc}(\R, L^2( \R^n) )$,
and $f$ is a nonlinear function
satisfying,  for every $ x \in \R^n$ and $  s \in \R$,
\be
\label{f1}
     f(x,s)s \le  - \alpha_{1} | s |^p + \phi_{1}(x)  \quad \mbox{for some } \
  p \ge 2,
\ee
\be
\label{f2}
    | f(x,s) | \le \alpha_{2} | s |^{p-1} + \phi_{2}(x),
\ee
\be
\label{f3}
\frac{\partial f}{\partial s}(x,s) \le  \alpha_3 ,
\ee
where  $\alpha_1, \ \alpha_2$ and    $ \alpha_3 $   are
all positive constants, $\phi_1 \in L^1(\R^n)$,
and $\phi_2 \in L^2(\R^n) \cap L^q (\R^n)$
with $  \frac{1}{p} + \frac{1}{q} = 1$.
Denote by
 $ F(x,s)= \int_{0}^{s} f(x, \tau) d \tau $. Then we
 assume that $F$ satisfies
\be
\label{F1}
-\phi_4(x) - \alpha_4 | s |^p \le F(x,s) \le -\alpha_5 | s |^p +\phi_3(x),
\ee
where $\alpha_4$ and $\alpha_5$ are positive constants
and $\phi_3, \phi_4 \in L^1(\R^n)$.

As in the case of bounded domains (see, e.g., \cite{tem1}),  it can be proved  that
 if   $g \in L^2_{loc} (\R, L^2 (\R^n))$
 and
\eqref{f1}-\eqref{F1}
 hold true,  then problem \eqref{rd1}-\eqref{rd2}   is well-posed
in  $\ltwo $, that  is, for   every  $\tau \in \R$
 and   $ u_\tau  \in  \ltwo
$,  there exists a unique  solution $  u  \in C( [\tau, \infty),  \ltwo
 ) \bigcap L^2(\tau, \tau +T; H^1(\R^n))
\bigcap L^p(\tau, \tau +T; L^p(\R^n))$
for every $T>0$. Further, the solution is continuous with respect
to   $u_\tau $ in  $ \ltwo  $.
To construct a cocycle  $\phi$ for problem  \eqref{rd1}-\eqref{rd2}, we
denote by
 $\Omega =\R$,  and   define a shift operator $\theta_t$
on $\Omega$
 for every
$t \in \R$
 by
$$
\theta_t  (\tau ) = t+ \tau, \quad \mbox{for all} \ \ \tau \in \R.
$$
Let $\phi$ be a mapping  from $\R^+ \times \Omega \times   \ltwo  $
to $  \ltwo$ given by
$$
\phi(t, \tau, u_\tau ) =
 u(t+\tau, \tau, u_\tau),
$$
 where  $ t \ge 0$,  $  \tau \in \R $,
$ u_\tau  \in \ltwo  $,  and
$ u $ is the solution of problem \eqref{rd1}-\eqref{rd2}.
By the uniqueness of solutions, we find that
for every $t, s \ge 0$, $\tau \in \R$ and
$ u_\tau  \in \ltwo  $,
$$
\phi (t+s, \tau,   u_\tau  ) =
\phi (t, s+ \tau,  (\phi(s, \tau,  u_\tau  ) ) ).
$$
Then  it follows  that $\phi$ is a continuous
$\theta$-cocycle  on $\ltwo $.
The purpose of this paper is to study
the existence of   pullback attractors
for $\phi$ in an appropriate phase space.

 Let  $E$ be a subset of   $ \ltwo $ and denote
   by
   $$ \| E \| = \sup\limits_{x \in E}
   \| x\|_{\ltwo  }.
   $$
 Suppose
   $D =\{ D(t) \}_{t\in \R}$  is  a family of
   subsets of $\ltwo  $  satisfying
 \be
 \label{basin_cond}
 \lim_{t \to  - \infty} e^{    \lambda  t} \| D( t) \|^2 =0,
 \ee
 where $\lambda$ is the positive constant
 appearing in \eqref{rd1}.
 Hereafter,  we  use ${\mathcal{D}}_\lambda$
 to denote
    the   collection of all
 families  of subsets of $L^2(\R^n)$ satisfying \eqref{basin_cond},
 that is,
 \be
 \label{D_lambda}
 {{\mathcal{D}}_\lambda = \{  D =\{ D(t) \}_{t\in \R}:
 D  \ \mbox{satisfies} \  \eqref{basin_cond} \} }.
\ee
 Throughout this paper, we assume the   following conditions for the external term:
  \be
 \label{gcond}
 \int_{-\infty}^\tau e^{\lambda \xi} \| g(\xi)\|^2  d \xi
 <  \infty, \quad \forall \  \tau  \in \R,
 \ee
 and
 \be
 \label{ginfinity}
\limsup_{k \to \infty}  \int_{-\infty}^\tau  \int_{|x| \ge k}  e^{\lambda \xi}
  |g(x, \xi) |^2 dx d\xi =0,
\quad \forall \  \tau  \in \R.
 \ee

 We remark that
 condition  \eqref{gcond} is useful   for  proving  existence of absorbing
 sets for problem \eqref{rd1}-\eqref{rd2}, while
the  asymptotically null condition
 \eqref{ginfinity} is  crucial  for establishing the asymptotic  compactness
 of solutions.
  Notice
   that conditions \eqref{gcond} and \eqref{ginfinity}
do not require that   $g$ be bounded in $L^2(\R^n)$
when $t \to \pm \infty$.  Particularly,
These assumptions do not  have  any restriction
 on  $g$ when $t \to +\infty$.

 It  follows  from   \eqref{ginfinity}   that
 for every $\tau \in \R$ and $\eta>0$, there is $K=K(\tau, \eta)>0$
 such that
\be
\label{ginfinity2}
 \int_{-\infty}^\tau  \int_{|x| \ge K}
e^{\lambda \xi}   |g(x, \xi) |^2 dx d\xi \le \eta e^{\lambda \tau} .
\ee
As we will see later,  inequality  \eqref{ginfinity2}  is
  crucial   for  deriving  uniform estimates on the tails of solutions
and these  estimates  are necessary  for
 proving  the asymptotic compactness of solutions.

\section{Uniform     estimates of solutions }
\setcounter{equation}{0}

      In this section, we
 derive uniform estimates of  solutions  of  problem \eqref{rd1}-\eqref{rd2}
 defined on $\R^n$
when $t \to \infty$.
 We  start with the estimates in $\ltwo  $.

\begin{lem}
\label{lem41}
 Suppose  \eqref{f1} and \eqref{gcond} hold.
Then for every $\tau \in \R$ and $D=\{D(t)\}_{t\in \R} \in {\mathcal{D}}_\lambda$,
 there exists  $T=T(\tau, D)>0$ such that for all $t \ge T$,
$$
\| u(\tau, \tau -t, u_0(\tau -t) ) \|^2  \le M + M  e^{- \lambda  \tau}
\int_{-\infty}^\tau
e^{\lambda \xi}
   \| g(\xi )\|^2    d \xi ,
$$
$$
\int_{\tau -t}^{\tau}  e^{\lambda  \xi } \|     u (\xi, \tau -t, u_0(\tau -t) )  \|^p_p
d\xi
  \le Me^{\lambda \tau} + M
    \int_{-\infty}^\tau
e^{\lambda \xi}
    \| g(\xi )\|^2    d \xi
 ,
 $$
and
 $$
\int_{\tau -t}^{\tau}  e^{\lambda  \xi } \|
 u (\xi, \tau -t, u_0(\tau -t) )  \|^2_{H^1}
d\xi
  \le Me^{\lambda \tau} + M
    \int_{-\infty}^\tau
e^{\lambda \xi}
    \| g(\xi )\|^2    d \xi
 ,
 $$
where
  $  u_0(\tau -t)  \in D(\tau -t)$,  and
$M$ is a  positive constant  independent of
$\tau$ and $D$.
\end{lem}

 \begin{proof}
 Taking  the inner product of  \eqref{rd1}
 with $u$ in $L^2(\R^n)$  we get that
 \be
 \label{p41_1}
 \frac{1}{2} \frac{d}{dt} \| u \|^2
+ \| \nabla u \|^2 + \lambda  \| u \|^2  = \int_{\R^n} f(x, u)  u dx + (g,u).
 \ee
 For the nonlinear term, by \eqref{f1} we have
\be
\label{p41_2}
 \int_{\R^n} f(x, u) u dx
\le  - \alpha_1 \int_{\R^n} |u|^p  dx + \int_{\R^n} \phi_1 dx.
\ee
By the Young inequality, the last term on the right-hand side
of \eqref{p41_1}  is bounded  by
\be
\label{p41_3}
|(g, u) | \le \| g \| \ \| u \|
\le {\frac 14} \lambda \|  u \|^2 + {\frac 1{ \lambda}} \| g \|^2.
\ee
It follows from \eqref{p41_1}-\eqref{p41_3} that
\be
\label{p41_4}
{\frac d{dt}} \| u \|^2 + 2 \| \nabla  u \|^2 + \lambda \| u \|^2
+  {\frac 12} \lambda \| u \|^2
+ 2 \alpha_1  \int_{\R^n} |u|^p dx \le C + {\frac 2\lambda} \| g \|^2.
\ee
Multiplying  \eqref{p41_4} by $e^{\lambda t}$ and then integrating
the resulting  inequality on $(  \tau -t, \tau)$ with $t \ge 0$, we find that
$$
\| u (\tau, \tau -t, u_0(\tau -t) ) \|^2
+ 2 e^{-\lambda \tau} \int_{\tau -t}^\tau e^{\lambda \xi} \| \nabla u (\xi, \tau -t, u_0(\tau -t )) \|^2 d\xi
$$
$$
+ {\frac 12} \lambda  e^{-\lambda \tau}
  \int^\tau _{\tau -t}  e^{\lambda \xi} \| u(\xi, \tau -t, u_0(\tau -t ) )\|^2 d\xi
+ 2 \alpha_1 e^{-\lambda \tau}
  \int^\tau _{\tau -t}  e^{\lambda \xi} \| u(\xi, \tau -t, u_0(\tau -t ) )\|^p_p d\xi
  $$
  $$
  \le e^{-\lambda \tau}  e^{\lambda  (\tau -t )}  \| u_0 (\tau -t ) \|^2
  + {\frac 2\lambda} e^{-\lambda \tau}
  \int^\tau_{\tau -t} e^{\lambda \xi} \| g(\xi ) \|^2 d\xi
  + {\frac C\lambda}
  $$
  \be
  \label{p41_5}
  \le e^{-\lambda \tau}  e^{\lambda  (\tau -t )}  \| u_0 (\tau -t ) \|^2
  + {\frac 2\lambda} e^{-\lambda \tau}
  \int^\tau_{-\infty} e^{\lambda \xi} \| g(\xi ) \|^2 d\xi
  + {\frac C\lambda}.
\ee
Notice that $u_0(\tau -t)  \in D(\tau -t)$
and $D= \{D(t)\}_{t \in \R} \in {\mathcal{D}}_\lambda$. We find that
for every $\tau \in \R$, there exists $T=T(\tau, D)$ such that
for all $t \ge T$,
\be
\label{p41_6}
e^{ \lambda (\tau-t) }   \| u_0(\tau -t ) \|^2
  \le
   {\frac {1}{   \lambda}}   \int_{-\infty}^\tau
e^{\lambda \xi}  \| g(\xi )\|^2 d \xi  .
\ee
By \eqref{p41_5}-\eqref{p41_6} we get that,
  for all $t \ge T$,
  $$
\| u (\tau, \tau -t, u_0(\tau -t) ) \|^2
+ 2 e^{-\lambda \tau} \int_{\tau -t}^\tau e^{\lambda \xi} \| \nabla u (\xi, \tau -t, u_0(\tau -t )) \|^2 d\xi
$$
$$
+ {\frac 12} \lambda  e^{-\lambda \tau}
  \int^\tau _{\tau -t}  e^{\lambda \xi} \| u(\xi, \tau -t, u_0(\tau -t ) )\|^2 d\xi
+ 2 \alpha_1 e^{-\lambda \tau}
  \int^\tau _{\tau -t}  e^{\lambda \xi} \| u(\xi, \tau -t, u_0(\tau -t ) )\|^p_p d\xi
  $$
  $$
  \le
  {\frac 3\lambda} e^{-\lambda \tau}
  \int^\tau_{-\infty} e^{\lambda \xi} \| g(\xi ) \|^2 d\xi
  + {\frac C\lambda},
  $$
which completes the proof.
  \end{proof}

  The following lemma  is  useful for deriving
  uniform estimates of solutions in $H^1(\R^n)$.

 \begin{lem}
\label{lem42}
 Suppose  \eqref{f1} and \eqref{gcond} hold.
Then for every $\tau \in \R$ and $D=\{D(t)\}_{t\in \R} \in {\mathcal{D}}_\lambda$,
 there exists  $T=T(\tau, D)>2$ such that for all $t \ge T$,
$$
\int_{\tau -2}^{\tau}  e^{\lambda  \xi } \|
   u (\xi, \tau -t, u_0(\tau -t) )  \|^2
d\xi
  \le Me^{\lambda \tau} + M
    \int_{-\infty}^\tau
e^{\lambda \xi}
    \| g(\xi )\|^2    d \xi
 ,
 $$
 $$
\int_{\tau -2}^{\tau}  e^{\lambda  \xi } \|   \nabla u (\xi, \tau -t, u_0(\tau -t) )  \|^2
d\xi
  \le Me^{\lambda \tau} + M
    \int_{-\infty}^\tau
e^{\lambda \xi}
    \| g(\xi )\|^2    d \xi
 ,
 $$
 and
 $$
\int_{\tau -2}^{\tau}  e^{\lambda  \xi } \|    u (\xi, \tau -t, u_0(\tau -t) )  \|^p_p
d\xi
  \le Me^{\lambda \tau} + M
    \int_{-\infty}^\tau
e^{\lambda \xi}
    \| g(\xi )\|^2    d \xi
 ,
 $$
where  $  u_0(\tau -t)  \in D(\tau -t)$,  and
$M$ is a  positive constant  independent of
$\tau$ and $D$.
\end{lem}

\begin{proof}
By \eqref{p41_4} we find that
$$
{\frac d{dt}} \| u \|^2
+ \lambda  \| u \|^2
\le C + {\frac 2\lambda} \| g \|^2.
$$
Let $s \in [\tau -2, \tau]$ and $t\ge 2$.
Multiplying  the above   by $e^{\lambda t}$ and integrating
over $(s, \tau -t)$,
  we  get
  $$
  e^{\lambda s} \| u(s, \tau -t, u_0(\tau -t) ) \|^2
  \le e^{\lambda (\tau -t)} \| u_0(\tau -t ) \|^2
  + C \int^s_{\tau -t} e^{\lambda \xi} d\xi
+ {\frac 2\lambda} \int^s_{\tau -t} e^{\lambda \xi} \| g(\xi)\|^2d\xi
$$
$$
 \le e^{\lambda (\tau -t)} \| u_0(\tau -t ) \|^2
  + {\frac C\lambda} e^{\lambda \tau}
+ {\frac 2\lambda} \int^\tau_{-\infty}
e^{\lambda \xi} \| g(\xi)\|^2d\xi.
$$
Therefore,
  there exists $T=T(\tau, D)>2$ such that
for all $t \ge T$ and $s \in [\tau -2, \tau]$,
\be
\label{p42_100}
  e^{\lambda s} \| u(s, \tau -t, u_0(\tau -t) ) \|^2
\le
{\frac C\lambda} e^{\lambda \tau}
+ {\frac 3\lambda} \int^\tau_{-\infty}
e^{\lambda \xi} \| g(\xi)\|^2d\xi.
\ee
Integrate the above with respect to $s$ on $(\tau -2, \tau)$
to obtain that
\be
\label{p42_101}
   \int_{\tau -2}^\tau
e^{\lambda s} \| u(s, \tau -t, u_0(\tau -t) ) \|^2 ds
\le
{\frac {2C}\lambda} e^{\lambda \tau}
+ {\frac 6\lambda} \int^\tau_{-\infty}
e^{\lambda \xi} \| g(\xi)\|^2d\xi.
\ee
On the other hand, for $s =\tau -2$, \eqref{p42_100}
implies  that
\be
\label{p42_1}
 e^{\lambda (\tau -2) } \| u(\tau -2, \tau -t, u_0(\tau -t) ) \|^2
\le
{\frac C\lambda} e^{\lambda \tau}
+ {\frac 3\lambda} \int^\tau_{-\infty}
e^{\lambda \xi} \| g(\xi)\|^2d\xi.
 \ee
Multiplying \eqref{p41_4} by $e^{\lambda t}$ and  then integrating
over $(\tau -2, \tau )$,  by \eqref{p42_1} we get that, for all $t\ge T$,
$$
e^{\lambda  \tau}  \| u(\tau, \tau -t, u_0(\tau -t) ) \|^2
 + 2    \int_{\tau -2}^\tau e^{\lambda \xi}
 \| \nabla u (\xi, \tau -t, u_0(\tau -t) ) \|^2 d\xi
 $$
 $$
  + 2  \alpha_1   \int_{\tau -2}^\tau e^{\lambda \xi}
 \|  u (\xi, \tau -t, u_0(\tau -t) ) \|^p_p  d\xi
 $$
 $$ \le
 e^{\lambda (\tau -2 )}    \| u(\tau -2, \tau -t, u_0(\tau -t) ) \|^2
 + {\frac {2}{\lambda}} \int_{\tau -2}^\tau e^{\lambda \xi}  \|  g (\xi) \|^2 d\xi
 + {\frac C\lambda} e^{\lambda \tau}
 $$
  $$
  \le C  \int_{-\infty}^\tau e^{\lambda  \xi}   \| g(\xi)\|^2  d\xi
  + C e^{\lambda \tau},
  $$
  which along with  \eqref{p42_101} completes the proof.
\end{proof}

 Note that $e^{\lambda \xi} \ge e^{\lambda \tau  -2\lambda} $
 for any $\xi \ge \tau -2$. So as an immediate  consequence
 of Lemma \ref{lem42} we have the following estimates.

 \begin{cor}
\label{cor43}
 Suppose  \eqref{f1} and \eqref{gcond} hold.
Then for every $\tau \in \R$ and $D=\{D(t)\}_{t\in \R} \in {\mathcal{D}}_\lambda$,
 there exists  $T=T(\tau, D)>2$ such that for all $t \ge T$,
 $$
\int_{\tau -2}^{\tau}    \|     u (\xi, \tau -t, u_0(\tau -t) )  \|^2
d\xi
  \le M  + M e^{- \lambda \tau}
    \int_{-\infty}^\tau
e^{\lambda \xi}
    \| g(\xi )\|^2    d \xi
 ,
 $$
$$
\int_{\tau -2}^{\tau}    \|   \nabla u (\xi, \tau -t, u_0(\tau -t) )  \|^2
d\xi
  \le M  + M e^{- \lambda \tau}
    \int_{-\infty}^\tau
e^{\lambda \xi}
    \| g(\xi )\|^2    d \xi
 ,
 $$
 and
 $$
\int_{\tau -2}^{\tau}    \|    u (\xi, \tau -t, u_0(\tau -t) )  \|^p_p
d\xi
  \le M  + M  e^{- \lambda \tau}
    \int_{-\infty}^\tau
e^{\lambda \xi}
    \| g(\xi )\|^2    d \xi
 ,
 $$
 where  $  u_0(\tau -t)  \in D(\tau -t)$,  and
$M$ is a  positive constant  independent of
$\tau$ and $D$.
\end{cor}

Next we derive uniform estimates of solutions
in $H^1(\R^n)$.

 \begin{lem}
\label{lem44}
 Suppose  \eqref{f1}, \eqref{F1}  and \eqref{gcond} hold.
Then for every $\tau \in \R$ and $D=\{D(t)\}_{t\in \R} \in {\mathcal{D}}_\lambda$,
 there exists  $T=T(\tau, D)>2$ such that for all $t \ge T$,
 $$
    \|   \nabla u (\tau, \tau -t, u_0(\tau -t) )  \|^2
  \le M + M e^{- \lambda \tau}
    \int_{-\infty}^\tau
e^{\lambda \xi}
    \| g(\xi )\|^2    d \xi  ,
 $$
 $$
  \|    u (\tau, \tau -t, u_0(\tau -t) )  \|^p_p
  \le M  + M e^{- \lambda \tau}
    \int_{-\infty}^\tau
e^{\lambda \xi}
    \| g(\xi )\|^2    d \xi  ,
 $$
 and
 $$
 \int^\tau_{\tau -1} \| u_\xi (\xi, \tau -t, u_0(\tau -t) ) \|^2 d\xi
 \le
 M  + M e^{- \lambda \tau}
    \int_{-\infty}^\tau
e^{\lambda \xi}
    \| g(\xi )\|^2    d \xi ,
 $$
 where  $  u_0(\tau -t)  \in D(\tau -t)$,  and
$M$ is a  positive constant  independent of
$\tau$ and $D$.
\end{lem}

\begin{proof}
In the following, we write $u_0(\tau -t)$ as $u_0$ for convenience.
  Taking the inner  product of \eqref{rd1}
  with $u_t$ in $L^2(\R^n)$ and then replacing $t$ by $\xi$,
we obtain
  $$
  \| u_{\xi} (\xi, \tau -t, u_0 ) \|^2 + {\frac d{d\xi}} \left (
  {\frac 12} \| \nabla u (\xi, \tau -t, u_0 )  \|^2
+ {\frac 12} \lambda \| u (\xi, \tau -t, u_0 )  \|^2
  -\int_{\R^n} F(x, u) dx \right )
  $$
  $$
  =(g(\xi), u_\xi (\xi, \tau -t,  u_0) ).
  $$
  Note that the  right-hand side of the above is bounded by
  $$
  |(g(\xi), u_\xi (\xi, \tau -t,  u_0) )|
  \le \| g(\xi) \|  \ \| u_\xi (\xi, \tau -t, u_0 ) \|
  \le
  {\frac 12} \| u_\xi (\xi, \tau -t, u_0 ) \|^2 +
  {\frac 12} \| g(\xi ) \|^2.
  $$
  Then we have
\be
\label{p44_1}
  \| u_{\xi} (\xi, \tau -t, u_0 ) \|^2 + {\frac d{d\xi}} \left (
   \| \nabla u (\xi, \tau -t, u_0 )  \|^2
+   \lambda \| u (\xi, \tau -t, u_0 )  \|^2
  - 2 \int_{\R^n} F(x, u) dx \right )
  \le \| g(\xi )\|^2,
\ee
which implies that
\be
\label{p44_2}
 {\frac d{d\xi}} \left (
   \| \nabla u (\xi, \tau -t, u_0 )  \|^2
+   \lambda \| u (\xi, \tau -t, u_0 )  \|^2
  - 2 \int_{\R^n} F(x, u) dx \right )
  \le \| g(\xi )\|^2.
\ee
Let $s \le \tau$  and $t\ge 2$.  By integrating
\eqref{p44_2} over $(s, \tau)$ we get
that
$$
 \| \nabla u (\tau, \tau -t, u_0 )  \|^2
+   \lambda \| u (\tau, \tau -t, u_0 )  \|^2
  - 2    \int_{\R^n} F(x, u (\tau, \tau -t, u_0 ) ) dx
  $$
 $$ \le
 \| \nabla u (s, \tau -t, u_0 )  \|^2
+   \lambda \| u (s, \tau -t, u_0 )  \|^2
  - 2    \int_{\R^n} F(x, u (s, \tau -t, u_0 ) ) dx
  + \int_s^\tau \| g(\xi) \|^2 d \xi.
  $$
  Now integrating the above with respect to $s$ on
  $(\tau -1, \tau)$ we  find that
  $$
 \| \nabla u (\tau, \tau -t, u_0 )  \|^2
+   \lambda \| u (\tau, \tau -t, u_0 )  \|^2
  - 2    \int_{\R^n} F(x, u (\tau, \tau -t, u_0 ) ) dx
  $$
  $$ \le
\int_{\tau -1}^\tau  \| \nabla u (s, \tau -t, u_0 )  \|^2  ds
+   \lambda \int_{\tau -1}^\tau\| u (s, \tau -t, u_0 )  \|^2 ds
 $$
\be
\label{p44_3}
 - 2  \int_{\tau -1}^\tau  \int_{\R^n}
 F(x, u (s, \tau -t, u_0 ) ) dx ds
  + \int_{\tau -1} ^\tau \| g(\xi) \|^2 d \xi.
\ee
By \eqref{F1} we  have
\be
\label{p44_4}
\alpha_5 \| u(\tau, \tau -t, u_0(\tau -t) )\|^p_p
-  \int_{\R^n} \phi_3 (x) dx
\le - \int_{\R^n} F(x, u (\tau , \tau -t, u_0 ) ) dx,
\ee
and
\be
\label{p44_5}
- \int_{\R^n} F(x, u (s , \tau -t, u_0 ) ) dx
\le
\alpha_4  \| u(s, \tau -t, u_0(\tau -t) )\|^p_p
+   \int_{\R^n} \phi_4 (x) dx.
\ee
It follows from \eqref{p44_3}-\eqref{p44_5}  that
 $$
 \| \nabla u (\tau, \tau -t, u_0 )  \|^2
+   \lambda \| u (\tau, \tau -t, u_0 )  \|^2
  + 2  \alpha_5     \| u (\tau, \tau -t, u_0 )  \|^p_p
  $$
  $$ \le
\int_{\tau -1}^\tau  \| \nabla u (s, \tau -t, u_0 )  \|^2  ds
+   \lambda \int_{\tau -1}^\tau\| u (s, \tau -t, u_0 )  \|^2 ds
 $$
 $$
  +  2  \alpha_4  \int^\tau_{\tau -1}
 \| u (s, \tau -t, u_0 (\tau -t)  )  \|^p_p ds
 + \int^\tau_{\tau -1} \| g(\xi) \|^2 d\xi
 + 2 \int_{\R^n}  (\phi_3 (x) + \phi_4(x) ) dx,
 $$
 which along with Corollary \ref{cor43} implies that there
 exists $T=T(\tau, D)>2$ such that for all
 $t \ge T$,
 $$
 \| \nabla u (\tau, \tau -t, u_0 )  \|^2
+   \lambda \| u (\tau, \tau -t, u_0 )  \|^2
  + 2  \alpha_5     \| u (\tau, \tau -t, u_0 )  \|^p_p
  $$
$$
  \le
  C+ C e^{-\lambda \tau} \int^\tau_{-\infty}
  e^{\lambda \xi} \| g(\xi) \|^2 d\xi
  + \int_{\tau -1}^\tau \| g(\xi) \|^2 d\xi
$$
 \be
 \label{p44_6}
  \le
  C+ C e^{-\lambda \tau} \int^\tau_{-\infty}
  e^{\lambda \xi} \| g(\xi) \|^2 d\xi
  +  e^\lambda e^{-\lambda \tau}
\int_{-\infty}^\tau e^{\lambda \xi}  \| g(\xi) \|^2 d\xi.
  \ee
  Similarly, first  integrating \eqref{p44_2}
  with respect to $\xi$ on $(s, \tau -1)$ and then
  integrating
  with respect to $s$ on $(\tau -2, \tau -1)$,
  by using Corollary \ref{cor43} we can get
  that for all $t \ge T$,
   $$
 \| \nabla u (\tau -1, \tau -t, u_0 )  \|^2
+   \lambda \| u (\tau-1, \tau -t, u_0 )  \|^2
  + 2  \alpha_5     \| u (\tau-1, \tau -t, u_0 )  \|^p_p
  $$
\be
\label{p44_7}
  \le
  C+ C e^{-\lambda \tau} \int^\tau_{-\infty}
  e^{\lambda \xi} \| g(\xi) \|^2 d\xi
  + \int_{\tau -2}^{\tau -1} \| g(\xi) \|^2 d\xi.
\ee
Now integrating \eqref{p44_1} over $(\tau -1, \tau)$ we obtain that
$$
\int^\tau_{\tau -1}
 \| u_{\xi} (\xi, \tau -t, u_0 ) \|^2 d\xi
 +   \| \nabla u (\tau, \tau -t, u_0 )  \|^2
+   \lambda \| u (\tau, \tau -t, u_0 )  \|^2
  - 2 \int_{\R^n} F(x, u (\tau) ) dx
  $$
  $$
  \le \int^\tau_{\tau -1}  \| g(\xi )\|^2 d\xi
  +   \| \nabla u (\tau-1, \tau -t, u_0 )  \|^2
+   \lambda \| u (\tau-1, \tau -t, u_0 )  \|^2
  - 2 \int_{\R^n} F(x, u (\tau-1) ) dx,
  $$
  which along with \eqref{p44_4}, \eqref{p44_5} and
  \eqref{p44_7} shows that for all $t \ge T$,
  $$
\int^\tau_{\tau -1}
 \| u_{\xi} (\xi, \tau -t, u_0 ) \|^2 d\xi
 \le \int^\tau_{\tau -1}  \| g(\xi )\|^2 d\xi
 + 2\int_{\R^n} (\phi_3(x) + \phi_4(x)) dx
 $$
 $$
  +   \| \nabla u (\tau-1, \tau -t, u_0 )  \|^2
+   \lambda \| u (\tau-1, \tau -t, u_0 )  \|^2
+ 2 \alpha_4  \| u (\tau-1, \tau -t, u_0 )  \|^p_p
 $$
$$
  \le
  C+ C e^{-\lambda \tau} \int^\tau_{-\infty}
  e^{\lambda \xi} \| g(\xi) \|^2 d\xi
  + \int_{\tau -2}^{\tau  } \| g(\xi) \|^2 d\xi
$$
\be
 \label{p44_10}
  \le
  C+ C e^{-\lambda \tau} \int^\tau_{-\infty}
  e^{\lambda \xi} \| g(\xi) \|^2 d\xi
+  e^{2\lambda} e^{-\lambda \tau}
\int_{-\infty}^\tau e^{\lambda \xi}  \| g(\xi) \|^2 d\xi.
\ee
Then Lemma \ref{lem44} follows  from \eqref{p44_6} and
\eqref{p44_10} immediately.
\end{proof}

We now  derive uniform estimates
of the derivatives of solutions in time. To this end, we also
assume ${\frac {dg}{dt}} \in L^2_{loc}(\R, L^2(\R^n) )$.

\begin{lem}
\label{lem45}
 Suppose  \eqref{f1}-\eqref{F1}
   and  \eqref{gcond}   hold.
   Let  ${\frac {dg}{dt}} \in L^2_{loc}(\R, L^2(\R^n) )$.
Then for every $\tau \in \R$ and $D=\{D(t)\}_{t\in \R} \in {\mathcal{D}}_\lambda$,
 there exists  $T=T(\tau, D)>2$ such that for all $t \ge T$,
 $$
  \| u_\tau (\tau, \tau -t, u_0(\tau -t) ) \|^2
 \le
 M  + M e^{- \lambda \tau}
    \int_{-\infty}^\tau
e^{\lambda \xi}
    \| g(\xi )\|^2    d \xi
+ M  \int^\tau_{\tau -1} \| g_\xi (\xi) \|^2 d\xi
 ,
 $$
 where  $  u_0(\tau -t)  \in D(\tau -t)$,  and
$M$ is a  positive constant  independent of
$\tau$ and $D$.
\end{lem}

\begin{proof}
 Let $u_t =v$ and differentiate  \eqref{rd1} with respect
 to $t$ to get that
 $$
 {\frac {\partial v}{\partial t}} - \Delta v
 +\lambda v
 ={\frac {\partial f}{\partial u}}(x,u)  v
 + g_t (x, t).
 $$
 Taking the inner product of the above with $v$ in $L^2(\R^n)$,
 we obtain
 \be
 \label{p45_1}
 {\frac 12} {\frac d{dt}} \| v \|^2
 + \| \nabla v \|^2 + \lambda \| v\|^2
 =\int_{\R^n}  {\frac {\partial f}{\partial u}}(x,u)  |v(x,t) |^2 dx
 + \int_{\R^n} g_t(x,t) v(x,t) dx.
 \ee
 By \eqref{f3} and the Young inequality, it follows from
 \eqref{p45_1} that
 \be
 \label{p45_2}
  {\frac d{dt}} \| v \|^2
  \le  2 \alpha_3  \| v \| ^2
 + {\frac 1\lambda} \| g_t(t) \|^2.
 \ee
 Let $s \in [\tau -1, \tau]$  and $t \ge 1$.    Integrating
 \eqref{p45_2} on $(s, \tau)$, by $v=u_t$ we get that
 $$
 \| u_\tau (\tau, \tau -t, u_0(\tau -t)) \|^2
 \le \| u_s(s, \tau -t, u_0(\tau -t) )\|^2
 $$
 $$
  + 2\alpha_3 \int_s^\tau \| u_\xi (\xi, \tau -t, u_0(\tau -t) ) \|^2 d\xi
  + {\frac 1\lambda} \int_s^\tau \| g_\xi (\xi) \|^2 d\xi.
  $$
  $$
    \le \| u_s(s, \tau -t, u_0(\tau -t) )\|^2
    + 2\alpha_3 \int_{\tau -1}^\tau \| u_\xi (\xi, \tau -t, u_0(\tau -t) ) \|^2 d\xi
  + {\frac 1\lambda} \int_{\tau -1}^\tau \| g_\xi (\xi) \|^2 d\xi.
  $$
  Now integrating the above with respect to $s$
  on $(\tau -1, \tau)$ we find that
   $$
 \| u_\tau (\tau, \tau -t, u_0(\tau -t)) \|^2
 \le  \int_{\tau -1}^\tau \| u_s(s, \tau -t, u_0(\tau -t) )\|^2 ds
 $$
 $$
 + 2\alpha_3 \int_{\tau -1}^\tau \| u_\xi (\xi, \tau -t, u_0(\tau -t) ) \|^2 d\xi
  + {\frac 1\lambda} \int_{\tau -1}^\tau \| g_\xi (\xi) \|^2 d\xi,
  $$
  which along with Lemma \ref{lem44} shows that there
  exists $T=(\tau, D)>2$ such that for all $t \ge T$,
   $$
 \| u_\tau (\tau, \tau -t, u_0(\tau -t)) \|^2
 $$
  $$
   \le
 C  + C e^{- \lambda \tau}
    \int_{-\infty}^\tau
e^{\lambda \xi}
    \| g(\xi )\|^2    d \xi
  + {\frac 1\lambda} \int_{\tau -1}^\tau \| g_\xi (\xi) \|^2 d\xi.
  $$
  The proof is completed.
\end{proof}

  We now  establish  uniform estimates on the tails
of solutions when $t \to \infty$. We show that the tails of solutions
are  uniformly  small  for  large space and time variables.
    These uniform estimates are crucial for proving
the pullback asymptotic compactness of the cocycle $\phi$.

 \begin{lem}
\label{lem46}
Suppose  \eqref{f1}, \eqref{F1}
   and  \eqref{gcond}-\eqref{ginfinity}   hold.
  Then for every $\eta>0$, $\tau \in \R$ and
 $D=\{D(t)\}_{t\in \R} \in {\mathcal{D}}_\lambda$,
 there exists  $T= T(\tau, D, \eta)>2$ and
$K=K(\tau, \eta)>0 $ such that for all $t \ge T$ and $k \ge K$,
$$
 \int _{|x| \ge  k}   |u(x, \tau, \tau -t, u_0(\tau -t) )|^2 dx
\le \eta,
  $$
where  $  u_0(\tau -t)  \in D(\tau -t)$,
  $K(\tau, \eta)$  depends   on $\tau$ and
  $\eta$,    and $T(\tau, D, \eta)$ depends
on  $\tau$, $D$ and   $\eta$.
\end{lem}

 \begin{proof}
  We  use a cut-off technique  to establish the estimates on the tails of solutions.
   Let $\theta$  be  a smooth function   satisfying   $0 \le \theta
(s) \le 1$ for $ s \in \R^+$, and
$$
\theta (s) =0 \ \mbox{for} \  0 \le s \le 1;
 \ \  \theta (s) =1 \
 \mbox{for} \ s \ge 2.
$$
Then there exists
  a constant
 $C$ such that  $ | \theta^{\prime}  (s) | \le C$
 for $  s \in \R^+ $.
Taking the inner product of \eqref{rd1} with
 $  \theta ({\frac {|x|^2}{k^2}}) u $ in  $\ltwo$, we get
$$
{\frac 12}    {\frac
 d{dt}} \int_{\R^n} \theta ({\frac {|x|^2}{k^2}}) |u|^2
  -     \int_{\R^n} \theta ({\frac {|x|^2}{k^2}}) u  \Delta u
  +    \lambda \int_{\R^n} \theta ({\frac {|x|^2}{k^2}}) |u|^2
$$
\be
\label{p466_1}
=  \int_{\R^n}  \theta ({\frac {|x|^2}{k^2}} )
  f(x,u)  u   dx
  +    \int_{\R^n} \theta ({\frac {|x|^2}{k^2}})  g(x,t) u(x,t) dx .
 \ee
 We now  estimate    the right-hand side of
\eqref{p466_1}.  For the nonlinear term,
by \eqref{f1}  we have
\be
\label{p466_2}
\int_{\R^n} \theta ({\frac {|x|^2}{k^2}} )
  f(x,u)  u   dx
  \le
  -\alpha_1 \int_{\R^n}  \theta ({\frac {|x|^2}{k^2}} ) |u|^p dx
  + \int_{\R^n}  \theta ({\frac {|x|^2}{k^2}} ) \phi_1 (x) dx
\le \int_{|x| \ge k}    \phi_1 (x) dx.
\ee
For the
  last  term on the right-hand side of \eqref{p466_1}
 we find that
 $$    \int_{\R^n}  \theta (\frac{|x|^2}{k^2})
     g(x, t)  u(x,t) dx
=   \int_{|x| \ge k} \theta (\frac{|x|^2}{k^2})
    g(x, t)  u(x,t) dx$$
  $$
\le
\frac12    \lambda  \int_{|x|\ge
k}\theta^2 (\frac{|x|^2}{k^2} )   |u  |^2 dx  +\frac
1{2\lambda}\int_{|x|\geq k} |g(x, t)|^2 dx $$
  \be
  \label{p466_3}
   \le
 \frac12     \lambda \int_{\R^n} \theta (\frac{|x|^2}{k^2})
  |u |^2  dx +\frac 1{2\lambda}
 \int_{|x|\geq k}  |g (x, t)
 |^2 dx .
 \ee
 On the other hand,
 for the  second  term on the  left-hand side of \eqref{p466_1},
 by integration by parts, we have
 $$
   \int_{\R^n} \theta ({\frac {|x|^2}{k^2}} ) u \Delta u
  = -   \int _{\R^n}  \theta
 ({\frac {|x|^2}{k^2}})  | \nabla u|^2
 -   \int_{\R^n}   \theta^{\prime} ({\frac {|x|^2}{k^2}})
 (  {\frac {2x}{k^2}} \cdot  \nabla u ) u .
 $$
 \be
 \label{p466_4}
\le
  -    \int _{ k \le |x| \le {\sqrt{2}}k}  \theta^{\prime} ({\frac {|x|^2}{k^2}})
 (  {\frac {2x}{k^2}} \cdot  \nabla u ) u
 \le  {\frac { C}k}  \int_{k \le |x| \le {\sqrt{2}}k } | u|  | \nabla u |
 \le
 {\frac {  C}{ k}} (  \| u \|^2 +  \| \nabla u \|^2 )
,
\ee
where $C $  is independent of    $k$.
It follows from  \eqref{p466_1}-\eqref{p466_4}
that
 $$
    {\frac d{dt}}
  \int_{\R^n} \theta ({\frac {|x|^2}{k^2}})
  |u|^2 dx
  +  \lambda    \int_{\R^n}  \theta ({\frac {|x|^2}{k^2}})
  |u|^2 dx
  $$
   \be
 \label{p466_5}
  \le   2  \int_{|x| \ge k} |\phi_1(x) | dx
  + {\frac 1{\lambda}} \int_{|x| \ge k} |g(x,t)|^2 dx
  + {\frac {C}{k}} ( \| u \|^2 + \| \nabla u \|^2 ).
\ee
Multiplying \eqref{p466_5} by $e^{\lambda t}$
and then
integrating over $(\tau -t, \tau)$ with $t\ge 0$, we get that
$$
 \int_{\R^n}  \theta ({\frac {|x|^2}{k^2}})
  |u(x, \tau, \tau -t, u_0(\tau -t) )|^2
   dx
  $$
  $$
  \le e^{-\lambda \tau}e^{\lambda (\tau - t)}
  \int_{\R^n}  \theta ({\frac {|x|^2}{k^2}})
 |u_0(x, \tau -t ) |^2   dx
  $$
  $$
  + 2  e^{-\lambda \tau} \int_{\tau -t}^\tau \int_{|x| \ge k}
e^{\lambda \xi}   |\phi_1 (x) | dx  d\xi
+ {\frac 1{\lambda}} e^{-\lambda \tau}
 \int_{\tau -t}^\tau \int_{|x| \ge k}
e^{\lambda \xi}   |g(x, \xi) |^2 dx d\xi
  $$
  $$
  + {\frac {C}{k}} e^{-\lambda \tau} \int_{\tau -t}^\tau \
e^{\lambda \xi}  \left (
 \| u (\xi, \tau -t, u_0(\tau -t ) )\|^2
 + \| \nabla u (\xi, \tau -t, u_0(\tau -t ) )\|^2
\right )
d \xi
$$
 $$
  \le e^{-\lambda \tau}e^{\lambda (\tau - t)}
  \int_{\R^n}
 |u_0(x, \tau -t ) |^2   dx
  $$
  $$
  + 2  e^{-\lambda \tau} \int_{-\infty}^\tau \int_{|x| \ge k}
e^{\lambda \xi}   |\phi_1 (x) | dx  d\xi
+ {\frac 1{\lambda}} e^{-\lambda \tau}
 \int_{-\infty}^\tau \int_{|x| \ge k}
e^{\lambda \xi}   |g(x, \xi) |^2 dx d\xi
  $$
\be
\label{p466_6}
  + {\frac {C}{k}} e^{-\lambda \tau} \int_{\tau -t}^\tau \
e^{\lambda \xi}  \left (
 \| u (\xi, \tau -t, u_0(\tau -t ) )\|^2
 + \| \nabla u (\xi, \tau -t, u_0(\tau -t ) )\|^2
\right )
d \xi.
\ee
Note that   given $\eta>0$, there is $T_1 =T_1(\tau, D, \eta)>0$
such that for all $t \ge T_1$,
\be
\label{p_466_7}
e^{-\lambda \tau}e^{\lambda (\tau - t)}
 \int_{|R^n}  |u_0( \tau -t )  |^2 dx   \le \eta.
  \ee
  Since $\phi_1 \in L^1(\R^n)$, there exists
  $K_1 =K_1(\eta)>0$ such that for all $k \ge K_1$,
  \be
  \label{p466_8}
  2  e^{-\lambda \tau} \int_{-\infty}^\tau \int_{|x| \ge k}
e^{\lambda \xi}   |\phi_1 (x) | dx  d\xi
\le \eta.
\ee
  On the other hand, by  \eqref{ginfinity2}
    there is $K_2=K_2(\tau, \eta)> K_1$ such that
  for all  $k \ge K_2$,
  \be
  \label{p_466_9}
  {\frac 1{\lambda}} e^{-\lambda \tau}
 \int_{-\infty}^\tau \int_{|x| \ge k}
e^{\lambda \xi}   |g(x, \xi) |^2 dx d\xi
\le  {\frac \eta\lambda}  .
\ee
For the last term on the right-hand side of
\eqref{p466_6},  it follows  from  Lemma \ref{lem41}
that there is
$T_2 =T_2(\tau, D)>0$ such that for all $t \ge T_2$,
$$
{\frac {C}{k}} e^{-\lambda \tau} \int_{\tau -t}^\tau \
e^{\lambda \xi}  \left (
 \| u (\xi, \tau -t, u_0(\tau -t ) )\|^2
 + \| \nabla u (\xi, \tau -t, u_0(\tau -t ) )\|^2
\right )
d \xi
$$
$$
\le {\frac {  C}{k}} \left (1+  e^{-\lambda \tau} \right )
\int_{-\infty}^\tau  e^{\lambda\xi}
  \| g(\xi) \|^2  d \xi.
$$
Therefore, there is  $K_3 =K_3(\tau, \eta)> K_2$ such that
for all $k\ge K_3$ and $t \ge T_2$,
\be
\label{p466_10}
{\frac {C}{k}} e^{-\lambda \tau} \int_{\tau -t}^\tau \
e^{\lambda\xi}  \left (
 \| u (\xi, \tau -t, u_0(\tau -t ) )\|^2
 + \| \nabla u (\xi, \tau -t, u_0(\tau -t ) )\|^2
\right )
d \xi \le \eta.
\ee
Let   $T=\max \{T_1, T_2 \}$.
Then by \eqref{p466_6}-\eqref{p466_10} we find that
for all $k\ge K_3$ and $t \ge T$,
$$
 \int_{\R^n}  \theta ({\frac {|x|^2}{k^2}})
  |u(x, \tau, \tau -t, u_0(\tau -t) )|^2 dx
  \le  3  \eta + {\frac \eta\lambda},
  $$
  and hence for all $k\ge K_3$ and $t \ge T$,
  $$
 \int _{|x| \ge \sqrt{2} k}
  |u(x, \tau, \tau -t, u_0(\tau -t) )|^2 dx
  $$
  $$ \le
 \int_{\R^n}  \theta ({\frac {|x|^2}{k^2}})
  |u(x, \tau, \tau -t, u_0(\tau -t) )|^2  dx
  \le    3  \eta + {\frac \eta\lambda},
  $$
  which completes the proof.
\end{proof}

 \section{Existence of pullback attractors}
\setcounter{equation}{0}

In this  section,  we prove
 the existence of a ${\mathcal{D}}_\lambda$-pullback
 global attractor
for the  non-autonomous   Reaction-Diffusion
equation  on $\R^n$.
We
 first  establish the   ${\mathcal{D}}_\lambda$-pullback asymptotic
compactness of  solutions   and  prove the existence
of a pullback attractor  in $L^2(\R^n)$.  Then we  show that
this attractor is actually a  ${\mathcal{D}}_\lambda$-pullback
attractor in $H^1(\R^n)$.

 \begin{lem}
\label{lem51}
Suppose  \eqref{f1}-\eqref{F1}
 and \eqref{gcond}-\eqref{ginfinity} hold.
Then $\phi$ is $\mathcal{D}_\lambda$-pullback
 asymptotically compact  in $\ltwo$,
that is, for  every $\tau \in \R$,
 $D=\{D(t)\}_{t\in \R} \in {\mathcal{D}}_\lambda$,
and $t_n \to \infty$,
 $u_{0,n}  \in D(\tau -t_n)$, the sequence
 $\phi(t_n, \tau -t_n,   u_{0,n}  ) $   has a
   convergent
subsequence in $\ltwo $.
\end{lem}

\begin{proof}
 We  use the uniform estimates on the tails of solutions to establish
 the precompactness of  $\phi  (t_n ,  \tau -t_n ,   u_{0,n}  )$
 in $\ltwo  $,  that is,  we   prove that for every
 $\eta>0$,  the sequence
 $\phi  (t_n ,  \tau -t_n ,   u_{0,n}  ) $ has a finite covering
 of balls of radii less than $\eta$.
 Given $K>0$,   denote by
\[{\Omega}_K = \{ x: |x|     \le K \} \quad \mbox{and} \quad
        {\Omega}^c_K = \{ x: |x|     >  K \}.\]
Then by  Lemma  \ref{lem46},    for the  given  $\eta >0$,
there exist $K= K(\tau, \eta)>0$  and $T=T(\tau, D, \eta)>2$
such that for $t \ge T$,
$$
\|   \phi  (t  ,  \tau -t  ,   u_{0}(\tau-t)   )  \| _{L^2 ({\Omega}^c_{K} )
   }
\le \frac{\eta}{4}.
$$
Since $t_n \to \infty$, there is $N_1=N_1(\tau, D, \eta)>0$ such that
$t_n \ge T$ for all $n \ge N_1$, and hence we obtain that,
for all  $n  \ge N_1$,
\be
\label{p51_1}
\|   \phi  (t_n  ,  \tau -t_n  ,   u_{0,n}    )  \| _{L^2 ({\Omega}^c_{K} )
 }
\le \frac{\eta}{4}.
\ee
On the other hand, by Lemmas  \ref{lem41} and \ref{lem44},    there
exist $C=C(\tau)>0$ and $N_2(\tau, D)>0$
such that for all $n \ge N_2 $,
\be
\label{p51_2}
\|   \phi  (t_n  ,  \tau -t_n  ,  u_{0,n}    )  \| _{H^1( {\Omega}_{K } )
} \le C.
\ee
By the compactness of embedding
$ H^1( {\Omega}_{K } )  \hookrightarrow  L^2 ( {\Omega}_{K } )$,
 the sequence $  \phi   (t_n  ,  \tau -t_n  ,  u_{0,n}    ) $ is precompact in
$L^2 ( {\Omega}_{K } )$.
Therefore,
for the given $\eta>0$,
$ \phi   (t_n  ,  \tau -t_n  ,   u_{0,n}   )$
has a finite covering in
$L^2 ( {\Omega}_{K } ) $ of
balls of radii less than $\eta/4$, which along with  \eqref{p51_1}
 shows
that $    \phi    (t_n  ,  \tau -t_n  ,  u_{0,n}    )   $ has a finite covering in
$ \ltwo   $ of balls of radii less than $\eta$, and
thus  $  \phi    (t_n  ,  \tau -t_n  ,   u_{0,n}    )$  is precompact
in   $\ltwo  $.
  \end{proof}

We now present  the    existence of a pullback   global
attractor
    for   $\phi$  in $L^2(\R^n)$.

 \begin{thm}
\label{thm52}
 Suppose  \eqref{f1}-\eqref{F1}
 and \eqref{gcond}-\eqref{ginfinity} hold.
Then problem \eqref{rd1}-\eqref{rd2}  has  a
 unique $\mathcal{D}_\lambda$-pullback  global attractor
$\{\mathcal{A}(\tau) \}_{\tau \in \R}\in {\mathcal{D}_\lambda} $  in $\ltwo
 $, that is,   for every  $\tau \in \R$,

(i) \  $\mathcal{A}(\tau)$ is compact in $L^2(\R^n)$;

(ii) \ $\{\mathcal{A}(\tau)\}_{\tau\in \R}$ is invariant, that is,
$$ \phi(t, \tau, \mathcal{A}(\tau)  )
= \mathcal{A}(t+\tau), \ \  \forall \   t \ge 0;
$$

(iii) \ \ $\{\mathcal{A}(\tau)\}_{\tau \in \R}$
attracts  every  set  in $\mathcal{D}_\lambda$  with respect to the norm
of  $L^2(\R^n)$,
 that is,  for every
 $B = \{B(\tau)\}_{\tau \in \R} \in \mathcal{D}_\lambda$,
$$ \lim_{t \to  \infty} d_{L^2(\R^n)}
 (\phi(t, {\tau -t}, B({\tau -t})), \mathcal{A}(\tau))=0,
$$
where   for any $Y, \ Z \subseteq  L^2(\R^n)$,
$$
d_{L^2(\R^n)} (Y,Z) =
  \sup_{y \in Y }
\inf_{z\in  Z}  \| y-z\|_{L^2(\R^n)}.
 $$
  \end{thm}

\begin{proof}
For $ \tau \in \R$, denote by
$$
B(\tau) = \{ u  : \ \| u \| ^2
\le  M+ M e^{-\lambda \tau} \int_{-\infty}^\tau
e^{\lambda \xi}   \| g(\xi) \|^2   d \xi \},
$$
where $M$ is the  positive constant in Lemma \ref{lem41}.
Note that $B=\{B(\tau)\}_{\tau \in \R} \in {\mathcal{D}_\lambda} $
is a ${\mathcal{D}_\lambda}$-pullback absorbing  for
$\phi$ in $\ltwo  $ by Lemma \ref{lem41}.
In addition, $\phi$ is ${\mathcal{D}_\lambda}$-pullback
asymptotically compact by Lemma \ref{lem51}. Thus the existence of
a ${\mathcal{D}_\lambda}$-pullback global attractor for $\phi$
in $L^2(\R^n)$  follows
from   Proposition \ref{att}.
  \end{proof}

In what follows, we strengthen Theorem \ref{thm52} and show that
the global attractor
$\{\mathcal{A}(\tau) \}_{\tau \in \R} $ is actually
a  $\mathcal{D}_\lambda$-pullback  global attractor
  in  $H^1(\R^n)$. As a necessary  step towards this goal,
  we  first prove the asymptotic compactness of solutions in
  $H^1(\R^n)$.

 \begin{lem}
\label{lem53}
Suppose  \eqref{f1}-\eqref{F1}
 and \eqref{gcond}-\eqref{ginfinity} hold. Let ${\frac {dg}{dt}} \in L^2_{loc} (\R, L^2(\R^n))$.
Then $\phi$ is $\mathcal{D}_\lambda$-pullback
 asymptotically compact in $H^1(\R^n)$,
that is, for  every $\tau \in \R$,
 $D=\{D(t)\}_{t\in \R} \in {\mathcal{D}}_\lambda$,
and $t_n \to \infty$,
 $u_{0,n}  \in D(\tau -t_n)$, the sequence
 $\phi(t_n, \tau -t_n,   u_{0,n}  ) $   has a
   convergent
subsequence in $H^1(\R^n)$.
\end{lem}

\begin{proof}
  By Lemma \ref{lem51},
the sequence
 $\phi(t_n, \tau -t_n,   u_{0,n}  ) =  u(\tau, \tau -t_n,   u_{0,n}  ) $ has a convergent
 subsequence in $L^2(\R^n)$, and hence
  there exists $v \in L^2(\R^n)$ such that,
  up to a subsequence,
  $$
  u(\tau, \tau -t_n,   u_{0,n}  )
  \to  v  \quad \mbox{in } \  \ L^2 (\R^n).
  $$
  This shows that
  \be
  \label{p53_1}
  u(\tau, \tau -t_n,   u_{0,n}  )
  \quad \mbox{is a Cauchy sequence  in }
  \ L^2(\R^n).
  \ee
  Next we prove that  $u(\tau, \tau -t_n,   u_{0,n}  )$ is actually
  a Cauchy sequence in $H^1(\R^n)$.
  For any $n, m \ge 1$,  it follows from \eqref{rd1} that
  $$
  -\Delta \left (
   u(\tau, \tau -t_n, u_{0,n})
   - u(\tau, \tau -t_m, u_{0,m})
  \right )
  + \lambda  \left (
   u(\tau, \tau -t_n, u_{0,n})
   - u(\tau, \tau -t_m, u_{0,m})
  \right )
  $$
  \be
  \label{p53_2}
 =
   f(x, u(\tau, \tau -t_n, u_{0,n}) )
   - f (x,  u(\tau, \tau -t_m, u_{0,m}) )
-
   u_\tau(\tau, \tau -t_n, u_{0,n})
   + u_\tau (\tau, \tau -t_m, u_{0,m}).
  \ee
 Multiplying \eqref{p53_2} by
 $u(\tau, \tau -t_n, u_{0,n})
   - u(\tau, \tau -t_m, u_{0,m})$,
   by \eqref{f3} we  get that
  $$
  \| \nabla  \left (
   u(\tau, \tau -t_n, u_{0,n})
   - u(\tau, \tau -t_m, u_{0,m})
  \right ) \|^2
  + \lambda  \|
   u(\tau, \tau -t_n, u_{0,n})
   - u(\tau, \tau -t_m, u_{0,m})
   \|^2
  $$
  $$
   \le \| u_\tau(\tau, \tau -t_n, u_{0,n})
   - u_\tau (\tau, \tau -t_m, u_{0,m}) \|
   \| u(\tau, \tau -t_n, u_{0,n})
   -u (\tau, \tau -t_m, u_{0,m})\|.
   $$
 \be
 \label{p53_3}
+
  \alpha_3
  \|
   u(\tau, \tau -t_n, u_{0,n})
   - u(\tau, \tau -t_m, u_{0,m})
   \|^2.
\ee
By Lemma \ref{lem45} we find that for every
$\tau \in \R$, there exists $T=T(\tau, D)$ such that
for all $t \ge T$,
$$
\| u_\tau(\tau, \tau -t , u_{0} (\tau -t) ) \|
\le  C.
$$
Since $t_n \to \infty$, there exists $N=N(\tau, D)$ such that
$t_n \ge T$ for all $n \ge T$.  Thus we obtain  that,
for all $n \ge N$,
$$
\| u_\tau(\tau, \tau -t_n , u_{0,n}  ) \|
\le  C,
$$
which along with \eqref{p53_3} shows that,
for all $n, m \ge N$,
$$
  \| \nabla  \left (
   u(\tau, \tau -t_n, u_{0,n})
   - u(\tau, \tau -t_m, u_{0,m})
  \right ) \|^2
  + \lambda  \|
   u(\tau, \tau -t_n, u_{0,n})
   - u(\tau, \tau -t_m, u_{0,m})
   \|^2
  $$
 $$
   \le  2C
   \| u(\tau, \tau -t_n, u_{0,n})
   -u (\tau, \tau -t_m, u_{0,m})\|.
   $$
 \be
 \label{p53_4}
+
  \alpha_3
  \|
   u(\tau, \tau -t_n, u_{0,n})
   - u(\tau, \tau -t_m, u_{0,m})
   \|^2.
\ee
It follows from \eqref{p53_1} and \eqref{p53_4}
that
$u(\tau, \tau -t_n, u_{0,n})$ is a Cauchy sequence
in $H^1(\R^n)$. The proof is completed.
\end{proof}

We are now ready to prove the existence of a global attractor
for  problem \eqref{rd1}-\eqref{rd2} in $H^1(\R^n)$.

 \begin{thm}
\label{thm54}
 Suppose  \eqref{f1}-\eqref{F1}
 and \eqref{gcond}-\eqref{ginfinity} hold.
 Let ${\frac {dg}{dt}} \in L^2_{loc} (\R, L^2(\R^n))$.
Then problem \eqref{rd1}-\eqref{rd2}  has  a
 unique $\mathcal{D}_\lambda$-pullback  global attractor
$\{\mathcal{A}(\tau) \}_{\tau \in \R}\in {\mathcal{D}_\lambda} $  in $H^1(\R^n)
 $, that is,   for every  $\tau \in \R$,

(i) \  $\mathcal{A}(\tau)$ is compact in $H^1(\R^n)$;

(ii) \ $\{\mathcal{A}(\tau)\}_{\tau\in \R}$ is invariant, that is,
$$ \phi(t, \tau, \mathcal{A}(\tau)  )
= \mathcal{A}(t+\tau), \ \  \forall \   t \ge 0;
$$

(iii) \ \ $\{\mathcal{A}(\tau)\}_{\tau \in \R}$
attracts  every  set  in $\mathcal{D}_\lambda$  with respect to the norm
of  $H^1 (\R^n)$,
 that is,  for every
 $B = \{B(\tau)\}_{\tau \in \R} \in \mathcal{D}_\lambda$,
$$ \lim_{t \to  \infty} d_{H^1 (\R^n)}
 (\phi(t, {\tau -t}, B({\tau -t})), \mathcal{A}(\tau))=0,
$$
where   for any $Y, \ Z \subseteq  H^1 (\R^n)$,
$$
d_{H^1 (\R^n)} (Y,Z) =
  \sup_{y \in Y }
\inf_{z\in  Z}  \| y-z\|_{H^1(\R^n)}.
 $$
  \end{thm}

  \begin{proof}
    The invariance of $\{\mathcal{A}(\tau)\}_{\tau \in \R}$
 is already given in Theorem \ref{thm52}.
So  we only need to prove (i) and  (iii).

   Proof of (i).  Let $\{v_n\}_{n=1}^\infty \subseteq \mathcal{A}(\tau)$. We want
   to show that  there exists $v \in \mathcal{A}(\tau)$ such
that, up to a subsequence,
 $ v_n   \to v $ in $H^1(\R^n)$.
 Since  $\mathcal{A}(\tau)$ is compact in $L^2(\R^n)$ by Theorem \ref{thm52},
 there exists $v \in \mathcal{A}(\tau)$ such that, up to a subsequence,
 \be
 \label{p54_1}
 v_n \to v \quad \mbox{in} \ \ L^2(\R^n).
 \ee
 We now prove the convergence  in \eqref{p54_1} actually
 holds in $H^1(\R^n)$.   Let $\{t_n\}_{n=1}^\infty$
be a sequence  with $t_n \to \infty$.
By the invariance of
 $\{\mathcal{A}(\tau)\}_{\tau \in \R}$, for every $n \ge 1$,
 there exists $w_n \in \mathcal{A} (\tau -t_n)$ such that
 \be
 \label{p54_2}
 v_n = \phi (t_n, \tau - t_n, w_n ).
\ee
By Lemma \ref{lem53}, it follows from \eqref{p54_2} that,
there exist $\tilde{v} \in H^1(\R^n)$ such that, up to a subsequence,
\be
 \label{p54_3}
 v_n = \phi (t_n, \tau - t_n, w_n ) \to \tilde{v} \quad \mbox{in } \ H^1(\R^n) .
\ee
Notice that \eqref{p54_1} and \eqref{p54_3} imply
$\tilde{v} =v \in \mathcal{A}(\tau)$, and thus (i) follows.

Proof  of (iii).  Suppose (iii) is  not true. Then there are
$\tau \in \R$,
$B = \{B(\tau)\}_{\tau \in \R} \in \mathcal{D}_\lambda$,
$\epsilon_0>0$ and $t_n \to \infty$ such that
$$
  d_{H^1 (\R^n)}
 (\phi(t_n, {\tau -t_n}, B({\tau -t_n})), \mathcal{A}(\tau)) \ge 2\epsilon_0,
$$
which implies that for every $n\ge 1$, there exists
$v_n \in B(\tau -t_n)$ such that
\be
\label{p54_10}
  d_{H^1 (\R^n)}
 (\phi(t_n, {\tau -t_n}, v_n ), \mathcal{A}(\tau)) \ge  \epsilon_0.
\ee
On the other hand,
By Lemma \ref{lem53}, there is $v \in H^1(\R^n)$ such that, up
to a subsequence,
\be
\label{p54_11}
\phi(t_n, {\tau -t_n}, v_n )
\to  v \quad \mbox{in } \ H^1(\R^n).
\ee
Since
$\{\mathcal{A}(\tau)\}_{\tau \in \R}$
attracts
 $B = \{B(\tau)\}_{\tau \in \R} $ in $L^2(\R^n)$ by Theorem \ref{thm52}, we have
\be
\label{p54_12}
 \lim_{n \to \infty } d_{L^2 (\R^n)}
 (\phi(t_n , {\tau -t_n }, v_n ), \mathcal{A}(\tau))=0.
\ee
By   \eqref{p54_11}-\eqref{p54_12}   and
 the compactness of $\mathcal{A}(\tau)$,
 we find that $v \in \mathcal{A}(\tau) $ and
 \be
\label{p54_15}
 \lim_{n \to \infty}  d_{H^1 (\R^n)}
 (\phi(t_n, {\tau -t_n}, v_n ), \mathcal{A}(\tau))
 \le \lim_{n \to \infty}  d_{H^1 (\R^n)}
 (\phi(t_n, {\tau -t_n}, v_n ),  v ) =0,
\ee
a contradiction with \eqref{p54_10}. The proof is completed.
  \end{proof}

 \end{document}